\documentclass{amsart}
\usepackage{amssymb} 
\usepackage{amsmath} 
\usepackage{amscd}
\usepackage{amsbsy}
\usepackage{comment}
\usepackage[matrix,arrow]{xy}
\usepackage{hyperref}

\DeclareMathOperator{\Norm}{Norm}

\DeclareMathOperator{\Rad}{Rad}
\DeclareMathOperator{\GL}{GL}

\DeclareMathOperator{\Gal}{Gal}

\DeclareMathOperator{\ord}{ord}

\newcommand{\Q}{{\mathbb Q}}

\newcommand{\F}{{\mathbb F}}

\def\mod#1{{\ifmmode\text{\rm\ (mod~$#1$)}
\else\discretionary{}{}{\hbox{ }}\rm(mod~$#1$)\fi}}

\begin {document}

\newtheorem{thm}{Theorem}
\newtheorem{lem}{Lemma}[section]

\newtheorem*{conj}{Conjecture}

\theoremstyle{definition}

\theoremstyle{remark}

\title[Erd\H{o}s--Selfridge]{Rational points on\\ Erd\H{o}s--Selfridge superelliptic curves}

\author[Michael Bennett]{Michael A. Bennett}
\address{Department of Mathematics, University of British Columbia, Vancouver, B.C., V6T 1Z2 Canada}
\email{bennett@math.ubc.ca}

\author{Samir Siksek}
\address{Mathematics Institute, University of Warwick, Coventry CV4 7AL, United Kingdom}
\email{S.Siksek@warwick.ac.uk}
\thanks{
The first author
is supported by NSERC, while the second author is supported by  
the EPSRC {\em LMF: L-Functions and Modular Forms} Programme Grant
EP/K034383/1, respectively.
}

\date{\today}


\begin {abstract}
Given $k \geq 2$, we show that there are at most finitely many rational numbers
$x$ and $y \neq 0$ and integers $\ell \geq 2$ 
(with $(k,\ell) \neq (2,2)$) for which
$$
x (x+1) \cdots (x+k-1) = y^\ell.
$$
In particular, if we assume that $\ell$ is prime,
 then all such triples $(x,y,\ell)$ satisfy either $y=0$ or $\log \ell < 3^k$.
\end {abstract}
\maketitle

\section{Introduction} \label{intro}

In a remarkable paper of 1975, Erd\H{o}s and Selfridge  \cite{ErSe} proved
that the product of at least two consecutive positive integers can never be a
perfect power. In other words, the Diophantine equation 
\begin{equation}
\label{main-eq}
x (x+1) \cdots (x+k-1) = y^\ell
\end{equation}
has no solutions in positive integers $x, y, k$ and $\ell$ 
with $k, \ell \geq 2$. Their proof, the culmination of more than forty years of
work by Erd\H{o}s, relied on an ingenious combination 
of elementary arguments and a lemma on bipartite
graphs.

For a fixed pair of positive integers $(k,\ell)$, 
equation \eqref{main-eq} defines a \emph{superelliptic curve} 
of genus at least $(\ell-1)(k-2)/2$.
In particular, if $\ell+k > 6$, the genus exceeds $1$, and by 
Faltings' theorem \cite{Fal}, the number of rational points $(x,y)$ is finite. 
Actually quantifying this result, for any given curve, can be an extremely
challenging problem. 

In the case of \emph{integer} points on superelliptic curves, one can typically
prove much stronger statements. In fact, 
given a polynomial $f(x)$ with integer
coefficients having at least two distinct roots,
a famous theorem of Schinzel and Tijdeman \cite{ST}
asserts that the integer solutions to the equation
$f(x)=y^\ell$ satisfy either $y \in \{ 0, \pm 1 \}$ or $\ell \leq \ell_0$ 
for some
(effectively computable) constant $\ell_0=\ell_0(f)$. 
Analogous absolute bounds upon exponents $\ell$ for which there exist non-trivial rational points on superelliptic curves are very hard to come by (though conjectured to exist). 
Indeed, such results for
the curves defined by equation \eqref{main-eq}, for small fixed values of $k$,
are among the very few in the literature (other results are restricted to polynomials of the shape $f(x)=g(h(x))$, where $g(x)=x^2+1$ or $x^3+1$ (see Darmon and Merel \cite{DM}) and to certain families of $g$ of small degree, treated in \cite{BD}).  These curves corresponding to \eqref{main-eq} admit a number of
obvious rational points, including ``trivial'' ones with $y=0$, and two
infinite families:
\begin{equation}  \label{ex1}
(x,y,k,\ell) = \left( \frac{a^2}{b^2-a^2}, \frac{ab}{b^2-a^2}, 2, 2 \right), \; \; a \neq \pm b
\end{equation}
and
\begin{equation}  \label{ex2}
(x,y,k,\ell) = \left( \frac{(1-2j)}{2} \, , \,  \frac{\pm 1}{2^j} \prod_{i=1}^j (2i-1) \, , \, 2j \, , \, 2 \right),
\end{equation}
where $a, b$ and $j$ are integers with $j$ positive. Two further solutions are given by
\begin{equation}  \label{ex3}
(x,y,k,\ell)=(-4/3,2/3,3,3) \quad \text{and} \quad  (-2/3,-2/3,3,3).
\end{equation}
It may be that there are no other such points and, in particular, none
whatsoever with $\ell \geq 4$. This is the content of a conjecture of Sander
\cite{Sa} (with requisite corrections noted in \cite{BBGH}) :

\begin{conj}[Sander] If $k \geq 2$ and $\ell \geq 2$ are integers, then the
only rational points on the superelliptic curve defined by equation
\eqref{main-eq} satisfy either $y=0$, 
or are as in \eqref{ex1}, \eqref{ex2} or \eqref{ex3}, for suitable choices
of the parameters $a, b$ and $j$.  
\end{conj}

Sander \cite{Sa} proved this conjecture for $2 \leq k \leq 4$ and, together
with Lakhal \cite{LaSa}, treated the case $k=5$. The conjecture was
subsequently established for $2 \leq k \leq 11$ by the first author, Bruin,
Gy\H{o}ry and Hajdu \cite{BBGH} (see also \cite{GHS}) 
and for $2 \leq k \leq 34$ by Gy\H{o}ry, Hajdu and Pint\'er \cite{GHP}.

In this short note, we will treat the case of arbitrary $k$. 
While we are not able to prove the above conjecture in its entirety, 
we establish the following partial result: 
\begin{thm} \label{main-thm}
Let $k \geq 2$ be a positive integer. 
Then equation \eqref{main-eq} has at most finitely many solutions in rational
numbers $x$ and $y$, and integers $\ell \geq 2$, with 
$(k,\ell) \neq (2,2)$ and $y \neq 0$. 
If we assume that $\ell$ is prime, all such solutions satisfy  
$\log \ell < 3^k$.  
\end{thm}

As far as the authors are aware, this is the first example of a rational analogue to the Schinzel-Tijdeman theorem to be proved for a superelliptic curve 
$f(x)=y^l$, where the polynomial $f$ has arbitrarily high degree and does not arise via composition from a polynomial of small degree.

\section{A Ternary Equation of Signature $(\ell,\ell,\ell)$}

\begin{lem}\label{lem:ternary}
Let $k \ge 2$ be an integer and $\ell > k$ be prime. Suppose
the superelliptic curve \eqref{main-eq} has an (affine)
rational point $(x,y)$ with $y \ne 0$.
Let $k/2<p \le k$ be prime.
Then there are non-zero integers $a$, $b$, $c$, $u$, $v$, $w$ satisfying
\begin{equation}\label{eqn:ternary}
a u^\ell + b v^\ell+ c w^\ell=0
\end{equation}
such that
\begin{enumerate}
\item[(i)] the integers $a$, $b$ and $c$ are $\ell$-th power free;
\item[(ii)] every prime divisor of $abc$ is at most $k$;
\item[(iii)] $p \nmid abc$; 
\item[(iv)] $p$ divides precisely one of $u$, $v$, $w$.
\end{enumerate}
\end{lem}
\begin{proof}
We write $x=n/s$ and $y=m/t$ where 
$m \ne 0$, the denominators $s$, $t$ are positive integers
and $\gcd(n,s)=\gcd(m,t)=1$. From equation \eqref{main-eq}, we have 
\[
\frac{n (n+s) (n+2s) \cdots (n+(k-1)s)}{s^k}=\frac{m^\ell}{t^\ell}.
\]
Our coprimality assumptions thus ensure that $s^k=t^\ell$. 
As $\ell$ and $k$ are coprime,
there is a positive integer $d$ such that $s=d^\ell$ and $t=d^k$. 
We are thus led to consider the equation 
\begin{equation} \label{main-eq2}
n (n+d^\ell) (n+2d^\ell) \cdots (n+(k-1)d^\ell) = m^\ell,
\end{equation}
where now all our variables are integers.
We write, for each $i \in \{ 0, 1, \ldots, k-1 \}$,
\begin{equation} \label{fool}
n + i d^\ell = a_i z_i^\ell,
\end{equation}
where $a_i$ is an $\ell$-th power free integer. 
Since the greatest common divisor of $n+i d^\ell$ and $n+j d^\ell$
divides $(i-j)$, each $a_i$ thus has the property
that its prime divisors  are bounded above by $k$.

Our argument relies on the basic fact that, given $k$ 
consecutive terms in arithmetic progression, each prime up to $k$ 
necessarily either  divides one of the terms, 
or the modulus of the progression. 
Fix a prime $p$ with $k/2 < p \leq k$. 

Suppose first that $p \mid d$. Then $p \nmid m$ and thus
$p \nmid a_i z_i$ for all $i$.
From \eqref{fool} we have 
\[
d^\ell+a_0 z_0^\ell-a_1 z_1^\ell=0;
\]
the proof of the lemma is complete in this case with
$a=1$, $b=a_0$, $c=-a_1$, $u=d$, $v=z_0$, $w=z_1$.

We may thus suppose $p \nmid d$. This fact combined with the inequality $p \le k$,
therefore forces $p$ to divide $n+id^\ell$
for some $0 \le i \le k-1$. Suppose first that $p$ does not
divide any other factor on the left-hand side of \eqref{main-eq2}.
Thus $p \nmid a_j z_j$ for $j \ne i$. Moreover,
$\ord_p(a_i z_i^\ell)=\ord_p(n+i d^\ell)=\ord_p(m^\ell)$ and so 
$p \nmid a_i$ and $p \mid z_i$ (as $a_i$ is $\ell$-th power free). 
By \eqref{fool} we have
\begin{gather*}
a_i z_i^\ell-a_{i+1} z_{i+1}^\ell+d^\ell=0 \qquad \text{if $i<k-1$}\\
a_{i} z_{i}^\ell-a_{i-1} z_{i-1}^\ell-d^\ell=0 \qquad \text{if $i=k-1$},
\end{gather*}
completing the proof in this case.

It remains to consider the case where $p$ divides at least two 
factors of the left-hand side of \eqref{main-eq2}.
In fact, as $p>k/2$ and $p \nmid d$,
precisely two factors are divisible by $p$ and these
have the form $n+i d^\ell$ and $n+(i+p) d^\ell$.
Thus $\ord_\ell((n+id^\ell)(n+(i+p)d^\ell)=\ord_\ell(m^\ell)$.
We shall make use of the identity
\[
(n+(i+p) d^\ell)(n+i d^\ell)-(n+(i+p-1) d^\ell)(n+(i+1) d^\ell)
+(p-1) d^{2\ell}=0.
\]
Substituting from \eqref{fool} completes the proof.
\end{proof} 

\section{Proof of Theorem~\ref{main-thm}}

We now turn to the proof of Theorem~\ref{main-thm}.
By previous work outlined in the introduction we may suppose that
$k \ge 35$. We shall suppose that $\ell>k$ is prime. Fix
a prime $k/2 <p \le k$ and suppose that \eqref{main-eq} has
a rational solution $(x,y)$ with $y \ne 0$. By 
Lemma~\ref{lem:ternary}, there are non-zero integers $a$, $b$, $c$,
$u$, $v$, $w$ satisfying \eqref{eqn:ternary} and conditions
(i)--(iv). By removing the greatest common factor, we may
suppose that the three terms in \eqref{eqn:ternary} are coprime
without affecting conditions (i)--(iv). After permuting the three
terms and changing signs if necessary, we may suppose further that
\[
a u^\ell \equiv -1 \pmod{4}, \qquad
b v^\ell \equiv 0 \pmod{2}.
\] 
Let $E$ be the Frey elliptic curve
\[
E \; : \; Y^2=X(X-a u^\ell) (X+b v^\ell).
\]
Write $G_Q=\Gal(\overline{\Q}/\Q)$. The action of 
$G_\Q$ on the $\ell$-torsion of $E$ gives rise to a
representation
\[
\overline{\rho}_{E,\ell} \; : \;
 G_\Q \rightarrow \GL_2(\F_\ell).
\]
As $\ell > k \ge 35$ and $E$ has full $2$-torsion,
we know by Mazur \cite{Mazur} 
that $\overline{\rho}_{E,\ell}$ is irreducible.
By the work of Kraus \cite{Kraus} (which appeals 
to modularity \cite{BreuilConradDiamondTaylor01} and Ribet's level lowering
\cite{Ribet-1990}) the representation $\overline{\rho}_{E,\ell}$ arises from
a weight $2$ newform $f$ of level $N^\prime$, where
\[
N^\prime=2^r \Rad_2(abc);
\]
here $r \le 5$ and $\Rad_2(n)$ denotes the product of the distinct
odd primes dividing $n$. By (ii) and (iii) of Lemma~\ref{lem:ternary}
we find that
\begin{equation} \label{enn}
N^\prime \mid  2^4 \cdot \prod_{q \leq k, q \neq p} q,
\end{equation}
where the product is over prime $q$.  
We appeal to the following standard result (see e.g. 
 \cite[Proposition 5.1]{Siksek}):
 \begin{lem}\label{lem:info}
Let $E/\Q$ be an elliptic curve of conductor $N$ and $f=q+\sum_{i \ge 2} c_i q^i$
be a newform of weight $2$ and level $N^\prime \mid N$.
Write $K=\Q(c_1,c_2,\dots)$ for the totally real number field
generated by the Fourier coefficients of $f$. 
If $\overline{\rho}_{E,\ell}$ arises from $f$
then there is some prime ideal $\lambda \mid \ell$
of $K$ such that for all primes $q$,
\begin{itemize}
\item if $q \nmid \ell N N^\prime$ then $a_q(E) \equiv c_q \pmod{\lambda}$;
\item if $q \nmid \ell N^\prime$ and $q \mid \mid N$ then
$q+1 \equiv \pm c_q \pmod{\lambda}$.
\end{itemize}
\end{lem}
Note that $\ell>k \ge p$ and so $\ell \ne p$.
Moreover, from \eqref{enn} we have $p \nmid N^\prime$.
Conclusion (iv) in Lemma~\ref{lem:ternary} ensures that $E$
has multiplicative reduction at $p$ and so $p \mid\mid N$.
We apply Lemma~\ref{lem:info} with $q=p$. Thus $\ell$ divides
$\Norm_{K/\Q}(p+1 \pm c_p)$. 
As $c_p$ (in any of the real embeddings
of $K$) is bounded by $2 \sqrt{p}$, this quantity is non-zero 
and hence provides an upper bound upon $\ell$:
\[
\ell \le (p+1+2\sqrt{p})^{[K:\Q]}=(\sqrt{p}+1)^{2[K:\Q]} .
\]
It remains to establish that $\log \ell < 3^k$.
The degree $[K:\Q]$
is bounded by $g_0^+(N^\prime)$
which denotes the dimension of the space of weight $2$ 
cuspidal newforms of level $N^\prime$.  
From Martin \cite{Ma}, we have, 
\[
g_0^{+}(N^\prime) \leq \frac{N^\prime+1}{12}.
\]
Thus
\[
\log{\ell}\le \frac{(N^\prime+1)}{6} \log{(\sqrt{p}+1)}. 
\]
By Schoenfeld \cite{Sc} 
$$
\sum_{\stackrel{q \leq k}{q \mbox{ \tiny{prime}}}} \log q < 1.000081 k.
$$
Finally, a routine computation making use of \eqref{enn}
and our assumption
$17 < k/2 \leq p \leq k$ allows us to conclude that $\log{\ell}< 3^k$.

\section{Concluding remark}

It is worth observing that our arguments employed to prove Theorem~\ref{main-thm} actually enable us to reach a like conclusion for curves of the shape
$$
x (x+1) \cdots (x+k-1) = b y^\ell,
$$
where $b$ is any integer with the property that its prime factors do not exceed $k/2$.




\end{document}